\documentclass{elsart}
\usepackage{graphicx,epsfig}
\usepackage{amsfonts}
\usepackage{amssymb}
\usepackage{amsmath}
\usepackage{rotating,subfigure}
\usepackage{multirow}
\usepackage{multicol,longtable,tabularx}

\begin{document}

\begin{frontmatter}

\title{Variable-step finite difference schemes for
the solution of Sturm-Liouville problems}
% \thanks[label1]{}
\author[label1]{Pierluigi Amodio\corauthref{cor1}},
\ead{pierluigi.amodio@uniba.it}
\author[label2]{Giuseppina Settanni},
\ead{giuseppina.settanni@unisalento.it}
\address[label1]{Dipartimento di Matematica, Universit\`a di Bari\\
via E. Orabona 4, 70125 Bari, Italy}
\corauth[cor1]{Corresponding author.}
\address[label2]{Dipartimento di Matematica e Fisica ``E. De Giorgi'', Universit\`a del Salento,
via per Arnesano, 73047 Monteroni di Lecce, Lecce, Italy}

\begin{abstract}
We discuss the solution of regular and singular Sturm-Liouville problems by means of High Order Finite Difference Schemes. We describe a code to define a discrete problem and its numerical solution by means of linear algebra techniques. Different test problems are proposed to emphasize the behaviour of the proposed algorithm.
\end{abstract}

\begin{keyword}
Second order BVPs, Sturm-Liouville problems, finite differences, high order, variable stepsize.
\MSC 65L15; 65L12; 65F15.
\end{keyword}

\end{frontmatter}

\section{Introduction}

Since many applications in quantum physics, quantum chemistry, science and industry are connected to Sturm-Liouville problems (SLPs), their solution has drummed up interest of several researches. Many codes have been developed to solve regular and singular problems. The first Fortran software package SLEIGN, introduced in \cite{Ba_th,BaGaKaZe,BaGoSh,Mar91}, computes automatically the eigenvalues and the eigenfunctions of some classes of SLPs, using the Pr\"ufer transformation and the oscillatory properties of the eigenfunctions. It requires the linear second order differential equation turned into a nonlinear first-order one of double size. Althought the code is reliable, sometimes it is slow and fails on some problems too. For this reason, a new version of this code, named SLEIGN2 \cite{BaEvZe}, has been developed in order to solve with any self-adjoint, separated or coupled boundary conditions, and with nearly all the classification of the endpoints of the interval, giving both quantitative and qualitative knowledge of the properties of self-adjoint Sturm-Liouville problems. Based on the approximation of the coefficients introduced by Pruess \cite{Pr}, two other Fortran codes have been also proposed. Both packages produce estimates of the eigenvalues and the eigenfunctions for regular and singular problems classifying the endpoints and the spectrum of the eigenvalues. In particular, the code SLEDGE in \cite{PrFu} fulfils the control of the global error, while SL02F in \cite{PrMa} uses a Pr\"ufer transformation.  Recently, it is implemented a MATLAB code MATSLISE, see \cite{LeDaBe}, applying Constant Perturbation Methods \cite{IxDmVb} for solving regular SLPs,  one-dimensional Schr\"odinger equations and radial Schr\"odinger equations with a distorted Coulomb potential. Actually, the first four codes have a larger range of applicability than the last one, although the latter is more efficient in the outlined class.

All these codes use standard ODE techniques to solve the discrete problem associated to SLPs. Across the last years, a class of matrix methods has been also developed, which applies finite differences or finite elements to SLPs, reducing the continuous problem to a matrix eigenvalue problem. In order to improve the numerical estimate of the eigenvalues, asymptotic corrections are also adopted, see \cite{AnPa, PaDeAn}.  Following this approach, in \cite{JAcGhGh,JAcGhMa,AcGhMa} a family of boundary value methods (BVMs) \cite {BrTr.book} has been introduced to approximate simultaneously the eigenvalues and the eigenfunctions of SLPs.

Based on this idea, in \cite{aip2010,jnaiam2011}  we illustrate a different approach able to solve regular and singular SLPs classifying the endpoints of the interval as regular or singular, limit points (LP) or limit circles (LC), oscillatory (LCO) or nonoscillatory (LCNO). The underlying methods approximate the derivatives of SLPs separately, therefore the second order differential equation is not transformed into an equivalent first order system; this means lower computational cost and easier stepsize selection strategy. The mainly used schemes are high order central finite differences, while other additional formulae are considered depending on the Sturm-Liouville classification of the endpoints of the interval.
%while an initial and final methods of the same order are applied on the extreme points of the interval. Obviously, the choose of the number of the initial and final methods depend on the Sturm-Liouville classification of the endpoints of the interval. The advantages of this approach are essentially two, the first is that no transformation of the original problem is required and the second is that to connect the solution of SLPs to a simple application of matrix methods for the eigenvalues.
The proposed approach with constant stepsize appears efficient and  robust to solve a particular singular Sturm-Liouville problem arising from the numerical computation of the eigenvalues and eigenfunctions of the finite (truncated) Hankel transform (see \cite{jamc2013}), important for numerous applications in signal and image process, and also for the solution of a two-parameter singular Sturm-Liouville problems derived from the numerical simulation of the so called `whispering gallery' modes (WGMs) occurring inside a prolate spheroidal cavity, see \cite{cpc2013}.

In this paper, stimulated by the previous results, we develop and test a code with variable stepsize and order, where the stepsize variation strategy takes advantage of an equidistribution of the error. Our aim is to point out the better efficency and the accuracy gained with the variable stepsize. For this reason in Section 2 we introduce the numerical schemes and give some hints on their usage to approximating SLPs. In Section 3 we sketch the code which is tested in the last section for some well known regular and singular problems.

\section{High Order Finite Difference Schemes}

In the last few years High Order Finite Difference Schemes have been largely used to solve different problems concerning with differential equations of order greater than 1.
We mention Initial Value Problems \cite{asc2012}, Boundary Value Problems \cite{jnaiam2009,aip2009,jcam2005,bit2007} and Sturm-Liouville problems \cite{jnaiam2011,jamc2013,cpc2013}. In all these cases the proposed formulae show good stability properties and the developed codes turn out to be competitive with the existing software, in particular when the problem to be solved is stiff and the solution has a slope different from that of its derivatives.

The reason of the nice behaviour depends on the good flexibility of the formulae which are chosen depending on the problem, the additional conditions and even the meshgrid and the point where the problem is discretized.

Suppose that the continuous solution $y(x)$ of the differential problem is defined on the interval $[a,b]$. Given a generic set of discretization points
\begin{equation} \label{grid}
X = \{ a = x_0 < x_1 < \dots < x_n = b \},
\end{equation}
we define the following $(s+r)$-steps formula in $x_i$
\begin{equation} \label{s-formule}
y^{(\nu)}(x_i) \simeq y_i^{(\nu)}
= \frac{1}{h_i^{\nu}} \sum_{j=-s}^r \alpha_{s+j}^{(\nu,s)} y_{i+j}, \qquad h_i = x_i-x_{i-1}
\end{equation}
which depends on the integers $s$ and $r$, that is, on the stencil $x_{i-s}, \dots, x_{i+r}$ associated to \eqref{s-formule}. The coefficients $\alpha_{s+j}^{(\nu,s)}$ are chosen such that the order $p$ of the formula is maximum, that is $p=s+r-\nu+1$ for general $s$ and $r$ and, due to the symmetry of \eqref{s-formule}, $p=s+r-\nu+2$ if $s=r$, $\nu$ is even and the stepsize $h_i$ is fixed over the stencil.

Starting from second order BVPs, given a general continuous problem $$ f(x,y,y',y'')=0, $$ we define  a discrete nonlinear system $$ f(x_i,y_i,y_i',y_i'')=0, \quad \text{for } i=1,\dots,n-1, $$ that, combined with the boundary conditions, allows to compute a unique solution. In compact form we could consider as discrete approximation of each derivative $Y^{(\nu)} = A_{\nu} \tilde Y$, where $Y$ and $Y^{(\nu)}$ contain the unknowns $y_i$ and the associated derivatives $y_i^{(\nu)}$, for $i=1,\dots,n-1$, $\tilde Y$ is the extension of $Y$ with the boundary values $y_0$ and $y_n$, and $A_{\nu}$ is the $(n-1) \times (n+1)$ matrix built with the coefficients of the approximations. Hence the continuous problem is approximated by its discrete counterpart $$f(X,Y,Y',Y'')=0.$$

Matrices $A_{\nu}$ have essentially a banded structure since the approximation in $x_i$ of each derivative makes use of values $y_j$ chosen around $y_i$. The best choice for $r$ and $s$ is $r=s$ even if in \cite{bit2007}, in order to improve stability properties for singular perturbation problems, it is suggested to approximate the first derivative with $s=r-2$ or $s=r+2$ depending on the sign of the term multiplying $y'$. We are forced to use a different strategy in the first and last gridpoints since we have not sufficient values to the left (or to the right) of $x_i$ to obtain a symmetric stencil.

Regular Sturm-Liouville problems
\begin{equation} \label{SLP}
\begin{cases}
p(x) y'' + q(x) y' + r(x) y = \lambda w(x) y \\
y(a) = y(b) = 0
\end{cases}
\end{equation}
are autonomous, their solution is zero in the boundary points of $X$ and contain an unknown parameter (eigenvalue) which must be determined in order nonnull solutions (eigenfunctions) exist. Since vectors $y \ne 0$ satisfying \eqref{SLP} are infinite, it is imposed that they satisfy a normalization condition
\begin{equation} \label{normal}
\int_a^b |y(x)|^2 dx = 1.
\end{equation}
All the eigenvalues are real. If they are defined in ascending order ($\lambda_0$ is the smallest one), then the eigenvector corresponding to $\lambda_k$ has $k$ intersections with the $x$-axis.

Based on the previous formulae, \eqref{SLP}-\eqref{normal} is approximated by the discrete eigenvalue problem
\begin{equation} \label{disc_SLP}
\begin{cases}
(P A_2 + Q A_1 + R) Y = \lambda W Y, \\
v^T Y = 1,
\end{cases}
\end{equation}
where $P$, $Q$, $R$ and $W$ are diagonal matrices, now $A_1$ and $A_2$ are square matrices (the first and the last columns are deleted because of the boundary conditions), and $v$ contains the coefficients of the quadrature formula discretizing \eqref{normal}. We remark that the order of this last formula is independent of the order of the overall method since it represents a scaling factor used to fix a particular eigenfunction.

If boundary conditions involve derivatives, then the extended vector $\tilde Y$ contains  also the values $y'_0$ and $y'_n$, requiring hence to compute new coefficients for the approximation. Taking advantage of the boundary conditions, the descrete problem is reformulated like in the regular case, see \cite{jnaiam2011}.

In the sequel we will also solve singular problems for which one or both the boundary conditions are useless. Supposing, for example, that $a$ is a LCNO point, we approximate the equation in \eqref{SLP} in $x_0$ with \eqref{s-formule} and $i=0$. Similarly, we proceed in the case $a$ is LP point and no boundary condition needs.

%In such a situation, we have to modify the previous strategy and use
%
%In the following we will choose even order and
%\begin{itemize}
%\item $r=s=p/2$ for the first derivative when $i=1,\dots,n$, ...
%\item $r=s-1$ or $r=s+1$ depending on the sign of the second derivative
%\end{itemize}
%
% For the computation of $\lambda$ we will use a proper algorithm for the eigenvalue computation.

\section{The HOFiD code for Sturm-Liouville problems}

The developed code for solving SLPs uses variable step/order to compute an approximation of a selected eigenvalue and its corresponding eigenfunction. In input, it requires the index $k$ of the eigenvalue $\lambda_k$ to be computed. As optional input, the orders of the schemes and the corresponding exit tolerances to be satisfied, and the length of the initial grid. The default values are 4, 6 and 8 as orders of convergence, and $1e-4$, $1e-6$ and $1e-8$, respectively as output tolerances for the solution. Error is approximated by using two different orders $p$ and $p+2$. We use relative error for the eigenvalue and absolute error for the eigenvector. The default number of initial points is $\max(20,5k)$.

The code, written in the MATLAB language, may be synthesized in the following few lines. It computes a first approximation with the order 2 and then iterates on the input orders. The solution obtained with order $p$ is used as starting approximation for the order $p+2$. For easy of understanding we have deleted all the controls on the input parameters and the computed solution.

\texttt{\hspace*{5mm} function [lam, y, hh] = solver( problem, k, tol, order, n0 ) \\
\hspace*{5mm} [lam, y] = init\_approx( problem, k, n0 ); \\
\hspace*{5mm} hh = (xn-x0)/n0*ones(n0-1,1); \\
\hspace*{5mm} for i = 1:length(order) \\
\hspace*{1cm} [lam, y, hh] = numer\_solut( problem, k, order(i), ... \\
\hspace*{15mm} tol(i), lam, y, hh ); \\
\hspace*{5mm} end}

It is based on the following main modules:
\begin{itemize}
\item \texttt{[lam,y]=init\_approx(problem,k,n)} computes an initial approximation for the $k$th eigenvalue, and the corresponding eigenfunction by using $n$ constant steps of the order 2 method. Since this matrix is not symmetric, the subroutine computes all the eigenvalues and then selects the $k$th checking the number of zeros of the corresponding eigenvector.
\item \texttt{[lam,y,hh] = numer\_solut(problem,k,ord,tol,hh,lam,y)} computes an approximate solution with a prescribed tolerance starting from a suitable meshgrid previously obtained together with the approximation of the eigenvalue and eigenvector. Order of the used method is fixed inside this subroutine.
\end{itemize}

This second subroutine is summarized as follows:

\texttt{\hspace*{5mm} function  [lam, y, hh] = numer\_solut( problem, k, order, ... \\
\hspace*{1cm} tol, hh, lam, y) \\
\hspace*{5mm} M = matrix( problem, order, hh ); \\
\hspace*{5mm} [lam0, y] = eig\_compute(M-lam*eye(n),y); \\
\hspace*{5mm} lam = lam+lam0; \\
\hspace*{5mm} e = err( problem, order, hh, lam, y ); \\
\hspace*{5mm} while norm(e)>tol \\
\hspace*{1cm} hh = equidistribute( hh, e ); \\
\hspace*{1cm} W = matrix( problem, order, hh ) \\
\hspace*{1cm} [lam0, u] = eig\_compute(M-lam I,y); \\
\hspace*{1cm} lam = lam+lam0; \\
\hspace*{1cm} e = err( problem, order, hh, lam, y ); \\
\hspace*{5mm} end}

It take advantages of the following subroutines
\begin{itemize}
\item \texttt{M = matrix(problem,order,hh)} computes the matrix $M$ of the problem discretized with fixed order and variable stepsize.
\item \texttt{[lam0,y] = eig\_compute(M,y)} computes a new approximation of the solution as the smallest eigenvalue of the matrix $M+\lambda I$ by means of the inverse power method.
\item \texttt{e = err(problem,order,hh,lam,y)} computes a discrete error function on the considered grid by using the subsequent even order to that used in \texttt{eig\_compute}.
\item \texttt{hh = equidistribute(hh,e)} computes a new grid by means of an equidistribution of the error function. In order to simplify the definition of the methods, only one stepsize variation is allowed in each stencil. This means that a variation in the stepsize is allowed every $r+s$ steps.
\end{itemize}

\section{Numerical experiments}

This section is devoted to illustrate the behaviour of the code on some regular and singular Sturm-Liouville problems, defined on bounded or unbounded intervals. All these tests show different difficulties that can be used to stress the code. In fact, for many other regular equations both the strategies of stepsize variation and order variation work very well.

When a problem is defined on an unbounded interval, then we perform a change of variable in order to obtain a new problem on a bounded one. Moreover,  as suggested in \cite{BaGoSh,Mar95,PrFu}, if the coefficients of the equation are not defined in the endpoints of the interval $[a,b]$, we consider the truncated interval $[\alpha, \beta]$ with $a< \alpha< \beta <b$.

In this section, depending on the problem and on the feature we want to test, we modify the default values introduced in the previous section. In particular we fix the order of the methods to compare the behaviour of the code on singular problems. In fact, it is known (see, for example, \cite{JAcGhMa}) that finite differences suffer of an order reduction on such equations.

When a method of order $p$ is used, we compute the relative error for the $k$th eigenvalue
\begin{equation}\label{relerr}
E_r(\lambda_k) = \frac{\left|\lambda_k^{(p)}- \lambda_k^{(p+2)}\right|}{\left|\lambda_k^{(p+2)}\right|},
\end{equation}
and the absolute error for the eigenfunction associated to $\lambda_k$
\begin{equation}\label{relabs}
E_a(y_k) = \left\| y_k^{(p)}- y_k^{(p+2)}\right\|_\infty.
\end{equation}

\paragraph*{Problem 1.}

The Mathieu equation \cite{Pry}
\begin{equation}
- y''(x) + c \cos(x) y(x) = \lambda  y(x), \qquad x \in [0,40],
\end{equation}
has regular boundary conditions $y(0)=y(40)=0$ and oscillatory coefficients. It is known that lower eigenvalues are grouped in clusters of 6, and more and tighter clusters appear as $c$ increases.

This problem is a numerical challenge for codes which have to compute one particular eigenvalue, and the information given by the eigenfunction is necessary to check if the obtained eigenvalue is the right one. Therefore, their estimate reveals the reliability and the efficiency of code described in Section 3.  Even if the problem is regular, we need a good starting accuracy to compute the first eigenvalue. Therefore, as shown in \tablename~\ref{Tab1} for $c = 5$, we apply  the order and stepsize variation strategy, starting with 251 equidistant mesh points.

For the first cluster of 6 eigenvalues we use orders 6, 8 and 10 to reach the exit tolerance $10^{-8}$. For the estimate of successive clusters we adopt the same order and stepsize variation strategy, but a different combination of orders and exit tolerances. Since the problem becomes easier, we start with order 4 (and then use orders 6 and 8) to reach the final exit tolerance $10^{-6}$, obtaining a good compromise between accuracy, number of mesh point and computational cost. As in general, the code reaches a better accuracy in computing the eigenvalues rather than the eigenfunctions. 

\begin{table}[bht]
\caption{\textit{Mathieu equation with $c=5$. Exit tolerances are $10^{-3}$, $10^{-6}$ and $10^{-8}$ for the orders 6, 8 and 10, respectively; $10^{-2}$, $10^{-4}$ and $10^{-6}$ for the orders 4, 6 and 8, respectively. $n_0=251$.}}
\begin{center}
\resizebox{0.9\columnwidth}{!}{%
\begin{tabular}{|cccrcc|}
\hline
%\multicolumn{6}{c}{$p = 6,8,10$,  \,  $tol = 10^{-3},10^{-6},10^{-8}$} \\
%\hline
%\hline
\multicolumn{1}{|c}{$k$}  &  \multicolumn{1}{c} {$p$} & \multicolumn{1}{c}{$\lambda_k$} & \multicolumn{1}{c}{$n$} & \multicolumn{1}{c}{$E_r(\lambda_k)$} & \multicolumn{1}{c|}{$E_a(y_k)$} \\
\hline
0 & 6-8-10 & -3.484238869351126e+00 & 981 & 1.08e-14 & 3.06e-09 \\
1 & 6-8-10 & -3.484221911373827e+00 & 1642 & 3.29e-14 & 7.46e-09\\
2 & 6-8-10 & -3.484197999007796e+00 & 986 & 1.15e-14 &  4.12e-09\\
3 & 6-8-10 & -3.484172609556845e+00 & 1315 &  1.33e-14 & 6.25e-09 \\
4 & 6-8-10 & -3.484151559702016e+00 & 1416 & 1.22e-13 & 9.90e-09\\
5 & 6-8-10 & -3.484139672740876e+00 & 1158 & 9.55e-14 & 9.08e-09\\
6 & 4-6-8 & -5.995435510621165e$-$01 & 695 & 9.94e-10 & 1.05e-07\\
12 & 4-6-8 & 1.932914885763969e+00 & 501 & 1.97e-09 & 1.67e-09\\
\hline
\end{tabular}}
%\caption{Mathieu equation with $c=5$. Exit tolerances are $10^{-3}$, $10^{-6}$ and $10^{-8}$ for the orders 6, 8 and 10, respectively; $10^{-2}$, $10^{-4}$ and $10^{-6}$ for the orders 4, 6 and 8, respectively. }
\end{center}\label{Tab1}
\end{table}

\begin{figure}[bht]
\begin{center}
\subfigure[Eigenfunction $y_0(x)$]
   {\includegraphics[width=0.45\textwidth]{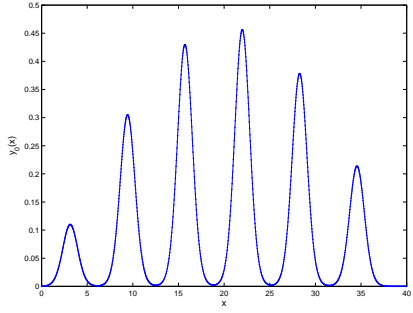}}
 \hspace{2mm}
 \subfigure[Absolute error and stepsize variation $y_0(x)$.]
   {\includegraphics[width=0.45\textwidth]{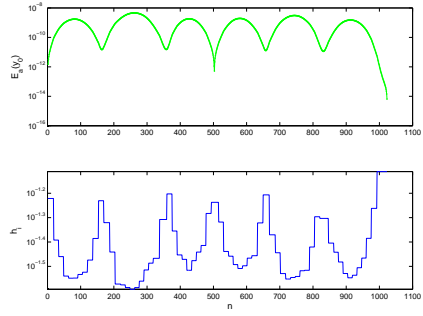}}
 \subfigure[Eigenfunction $y_5(x)$]
   {\includegraphics[width=0.45\textwidth]{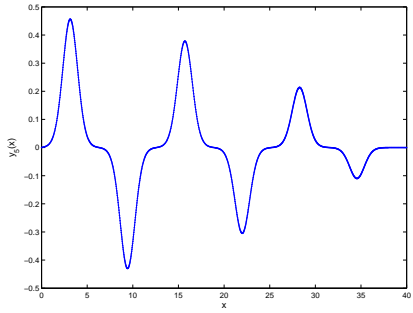}}
\hspace{2mm}
 \subfigure[Absolute error and stepsize variation $y_5(x)$.]
   {\includegraphics[width=0.45\textwidth]{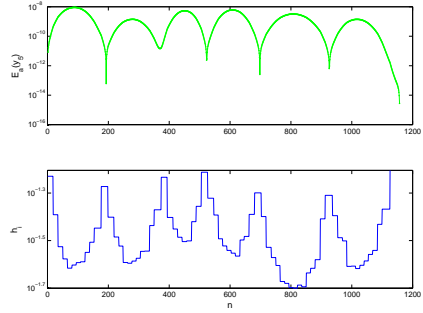}}
\caption{\textit{Mathieu equation with $c=5$. Solution, absolute error and stepsize variation for the first and sixth eigenfunction.}}\label{Fig1}
\end{center}
\end{figure}

In  \figurename~\ref{Fig1} we depict the eigenfunctions associated to the first and the last eigenvalue belonging to the first cluster, their absolute error and the stepsize variation. We observe how in this case it is quite important to use small exit tolerances since the numerical solution is often quite near the $x$-axis.

\paragraph*{Problem 2.}

The Pruess equation
\begin{equation}
- y''(x) + \mathrm{ln}(x) y(x) = \lambda  y(x), \qquad x \in [0,4],
\end{equation}
has regular boundary conditions $y(0)=y(4)=0$. It is defined as a \textit{regular problem that looks singular} in \cite{Pry} since for any high order finite difference scheme the order of convergence reduces to 2 as when these formulae are applied to singular problems. Consequently, it is simple to infer that there is no advantage in their use. 

In \tablename~\ref{Tab2} we just analyze the convergence behaviour for fixed order. We depict the required number of steps for the first, the fifth and the tenth eigenvalue and eigenfunction using orders 4, 6 and 8 and stepsize variation with the exit tolerance $10^{-8}$. As the results highlight, the convergence is not safe from the order reduction. Actually, order 4 is sufficient to guarantee a good accuracy with a small number of mesh points for the first eigenvalue, but when the index of the eigenvalue increases, then order 8 becomes competitive. Obviously, for simple problems the use of one order is enough to gain the accuracy of the lowest eigenvalues, while an order variation is recommended for the larger ones.

\begin{table}[bht]
\caption{\centerline{\textit{Pruess equation. Exit tolerance $10^{-8}$. $n_0=21$.}}}
\label{Tab2}
\begin{center}
\resizebox{0.8\columnwidth}{!}{%
\begin{tabular}{|cccrcc|}
\hline
$k$ & $p$ & $\lambda_k$ & $n$ & $E_r(\lambda_k)$ & $E_a(y_k)$\\
\hline
1 & 4 & 1.12481680 & 290 & 3.36e-09 & 2.64e-09\\
  & 6 & 1.12481680 & 435 & 1.44e-09 & 1.74e-09\\
  & 8 & 1.12481678 & 321 & 4.13e-09 & 5.18e-09\\
\hline
4 & 4 & 15.8644571 & 1201 & 1.13e-09 & 1.39e-09\\
  & 6 & 15.8644571 & 363 & 1.47e-09 & 7.65e-09 \\
  & 8 & 15.8644571 & 395 & 9.41e-10 & 5.38e-09 \\
\hline
9 & 4 & 62.0987975 & 1930 & 3.22e-10 & 6.74e-09 \\
  & 6 & 62.0987973 & 503 & 8.99e-10 & 3.91e-09 \\
  & 8 & 62.0987972 & 429 & 7.68e-10 & 8.85e-09\\
\hline
\end{tabular}}
\end{center}
\end{table}

\paragraph*{Problem 3.}

The Airy equation \cite{Pry}
\begin{equation}\label{Prob3}
- y''(x) + x y(x) = \lambda  y(x), \qquad x \in [0,+\infty],
\end{equation}
is a singular problem with boundary condition $y(0)=0$ while $b=+\infty$ is LP. The  eigenvalues are the zeros of the Airy function $$Ai(\lambda) = (J_{1/3}+J_{-1/3}) \left(\displaystyle{\frac{2}{3} \lambda^{1/3} }\right),$$ where $J_{\alpha}$ is the Bessel function. Since the problem is defined on a semi-infinite interval, we perform the following transformation 
\begin{equation}\label{var}
\tau= 1 - \frac{1}{\sqrt{1+x}} \in [0,1].
\end{equation}
Consequently, with the change of variable $u(\tau) = y(\tau(x))$, problem \eqref{Prob3} is rewritten as
\begin{equation}\label{Prob3u}
- \displaystyle{\frac{(1-t)^6}{4}} u''(\tau) + \frac{3}{4} (1-\tau)^5 u'(\tau) + \displaystyle{\frac{2 t -t^2}{1-t^2}} u(\tau) =  \lambda  u(\tau), \qquad \tau \in [0,1],
\end{equation}
regular in $a=0$, while $b=1$ is LP. Moreover, since the problem is singular in 1, the interval is truncated to $\beta = 1-\delta$ with $\delta=1e-4$.
\begin{table}[bht]
\caption{\centerline{\textit{Airy equation. Exit tolerance $10^{-8}$. $n_0=21$.}}}
\label{Tab3}
\begin{center}
\resizebox{0.7\columnwidth}{!}{%
\begin{tabular}{|cccccc|}
\hline
%\multicolumn{5}{|c|}{$k = 0$}\\
%\hline
$k$ & $p$ & $\lambda_k$ & $n$ & $E_r(\lambda_k)$ & $E_a(y_k)$\\
\hline
0 & 4 & 2.33810740 & 872 & 1.00e-09 & 4.98e-09\\
  & 6 & 2.33810741 & 789 & 1.36e-09 & 8.51e-09\\
  & 8 & 2.33810741 & 912 & 1.60e-09 & 1.00e-08\\
  & 4-6-8 & 2.33810741 & 938 & 1.19e-09 & 7.44e-09\\
\hline
%\multicolumn{5}{|c|}{$k = 4$}\\
%\hline
4 & 4 & 7.94413358 & 3648 & 2.00e-11 & 4.54e-09\\
  & 6 & 7.94413358 & 3657 & 5.28e-11 & 1.79e-09\\
  & 8 & 7.94413358 & 5169 & 4.33e-11 & 1.47e-09\\
  & 4-6-8 & 7.94413358 & 5337 & 7.94e-11 & 5.99e-11\\
            \hline
\end{tabular}}
\end{center}
\end{table}

As illustrated in \tablename~\ref{Tab3}, it is clear that in this example we have a reduction of the order of the convergence. Effectively order 4 is the best choice both for accuracy and minimum computational cost. This behavior justifies the failure of the order variation strategy, and underlines as it makes sense to apply only a variable stepsize stategy with the smaller order.

\paragraph*{Problem 4.}

The Laguerre's equation \cite{Pry}
\begin{equation} \label{pry_test}
- y''(x) + \left( x^2+\displaystyle{\frac{3}{4 x^2}} \right)  y(x) = \lambda  y(x), \qquad x \in [0,+\infty],
\end{equation}
is singular and both the endpoints $a = 0$ and $b = +\infty$ are LP. The eigenvalues satisfy the relation $\lambda_k = 4 (k+1)$, for $k\ge 0$. As in the previous problem, the equation is defined on a semi-infinite interval, therefore a change of variable is need to reformulate \eqref{pry_test} in the finite interval $[0,1]$. Moreover, due to the singularities in 0 and 1, the new problem is solved in the truncated interval $[\delta, 1-\delta]$, with $\delta=10^{-4}$.

As shown in \tablename~\ref{Tab4}, the results for this problem are completely different from the previous ones since now order is preserved and order (and stepsize) variation allow to greatly reduce the required number of points.

\begin{table}[bht]
\caption{\centerline{\textit{Laguerre's equation. Exit tolerance $10^{-8}$. $n_0=21$.}}}
\label{Tab4}
\begin{center}
\resizebox{0.7\columnwidth}{!}{%
\begin{tabular}{|cccccc|}
\hline
$k$& $p$ & $\lambda_k$ & n & $E_r(\lambda_k)$ & $E_a(y_k)$\\
\hline
%\multirow{5}{*}{0} & & & & & \\
0& 4 & 3.99999999  & 1018 & 3.65e-10 & 8.71e-09\\
& 6 & 4.00000000  & 969  & 1.19e-11 & 8.63e-09\\
&8 & 3.99999999 & 748 & 3.03e-11 & 9.04e-09\\
%&10& 3.99999999 & 763 & 3.50e-11 & 8.89e-09\\
&4-6-8 & 3.99999999 & 569 & 2.22e-11 & 6.77e-09\\
\hline
%\hline
4 &4-6-8 & 20.00000000 & 597& 1.33e-11 & 9.04e-09\\
\hline
%\hline
9 &4-6-8 & 39.99999999 & 753 & 2.95e-11 & 9.37e-09\\
\hline
%\hline
24 &4-6-8 & 99.99999999 & 2067 & 8.52e-12 & 1.09e-09\\
 \hline
\end{tabular}}
\end{center}
\end{table}

\begin{figure}[bht]
\begin{center}
\subfigure[Eigenfunction $u_9(\tau)$ in the interval $(0,1)$.]
   {\includegraphics[width=0.8\textwidth]{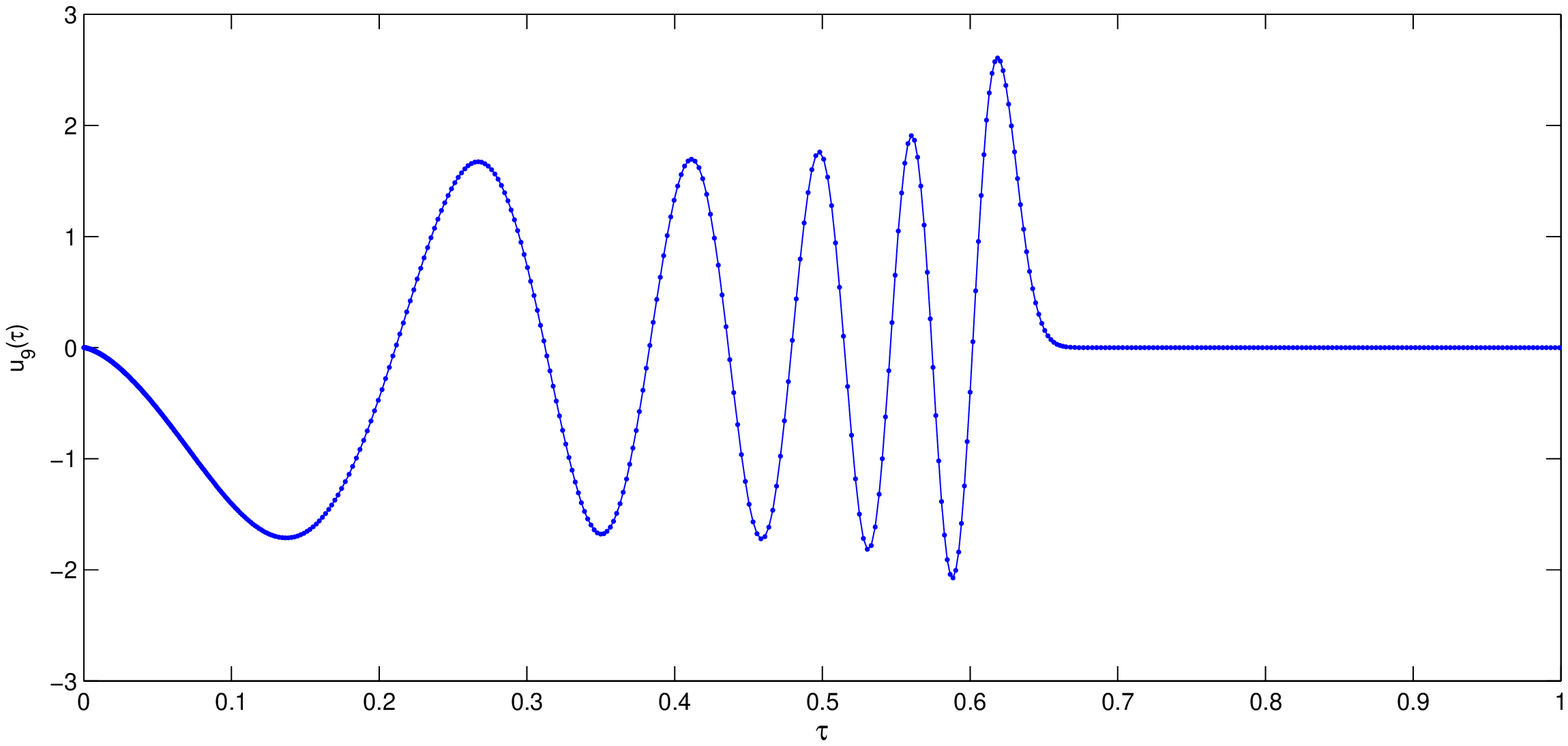}}
 \subfigure[Eigenfunction $y_9(x)$ in the truncated interval $(0,20)$.]
   {\includegraphics[width=0.8\textwidth]{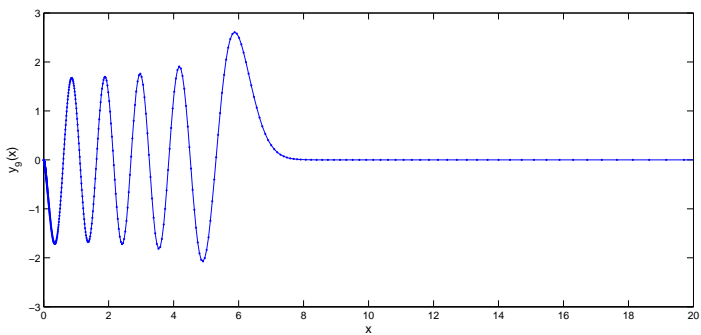}}
\end{center}
\caption{\textit{Laguerre's equation. Numerical approximation of the tenth eigenfunction.}}\label{Fig2}
\end{figure}

In \figurename~\ref{Fig2} we depict the tenth eigenfunction, both for the transformed problem and the original one truncated in the interval $ [0,20]$, where it is clear that stepsize variation is necessary to gain high accuracy with a limited number of steps.

\section{Conclusions}
In this paper we propose a code for solving regular and singular Sturm-Liouville problems which takes advantage of stepsize and order variation strategies. The code works well for any considered problem even if the use of high orders is sometimes useless because of a known deficiency of finite differences applied to singular equations. The successive step of this paper will be the application of this code to some challenging problems for which the use of a robust stepsize variation is necessary.


\begin{thebibliography}{99}

\bibitem{JAcGhGh}
L.~Aceto, P.~Ghelardoni and G.~Gheri,
\emph{An algebraic procedure for the spectral corrections using the miss-distance functions in regular and singular Sturm-Liouville problems},
SIAM J.\ Numer.\ Anal. \textbf{44} (2006), 2227--2243.


\bibitem{JAcGhMa}
L.~Aceto, P.~Ghelardoni and C.~Magherini,
\emph{BVMs for Sturm-Liouville eigenvalue estimates with general boundary conditions}, %
J.\ Numer.\ Anal.\ Ind.\ Appl.\ Math.\ \textbf{4} (2009), 113--127.

\bibitem{AcGhMa}
L.~Aceto, P.~Ghelardoni and C.~Magherini,
\emph{Boundary value methods as an extension of Numerov's method for Sturm-Liouville eigenvalue estimates},
Appl.\ Numer.\ Math.\ \textbf{59} (7) (2009), 1644--1656.

\bibitem{asc2012}
P.~Amodio, C.~Budd, O.~Koch, G.~Settanni and E.B.~Weinm\"uller,
\emph{Computations for a Model of Flow in Concrete},
ASC Report \textbf{26/2012}, Institute for Analysis and Scientific Computing, Vienna University of Technology, Wien.

\bibitem{jnaiam2009}
P.~Amodio and G.~Settanni, %
\emph{Variable step/order generalized upwind methods for the numerical solution of second order singular perturbation problems}, %
J.\ Numer.\ Anal.\ Ind.\ Appl.\ Math.\ \textbf{4} (2009), 65--76.

\bibitem{aip2009}
P.~Amodio and G.~Settanni, %
\emph{A deferred correction approach to the solution of singularly perturbed BVPs by high order upwind methods: implementation details}, %
in: Numerical analysis and applied mathematics - ICNAAM 2009. T.E. Simos, G. Psihoyios and Ch. Tsitouras (eds.), AIP Conf.\ Proc.\ \textbf{1168}, issue 1 (2009),
711--714.

\bibitem{aip2010}
P.~Amodio and G.~Settanni, %
\emph{High order finite difference schemes for the numerical solution of eigenvalue problems for IVPs in ODEs}, %
in: Numerical analysis and applied mathematics - ICNAAM 2010. T.E. Simos, G. Psihoyios and Ch. Tsitouras (eds.), AIP Conf.\ Proc.\ \textbf{1281}, issue 1 (2010), 202--204.

\bibitem{jnaiam2011}
P.~Amodio and G.~Settanni,
\emph{A matrix method for the solution of Sturm-Liouville problems},
JNAIAM J.\ Numer.\ Anal.\ Ind.\ Appl.\ Math.\ \textbf{6} (2011), 1--13.

%\bibitem{aip2011}
%P.~Amodio and G.~Settanni, \emph{A stepsize variation strategy for the solution of regular Sturm-Liouville problems},
%in: Numerical Analysis and Applied Mathematics - ICNAAM 2011. T. E. Simos, G. Psihoyios and Ch. Tsitouras editors,
%AIP Conf.\ Proc.\ \textbf{1389}, issue B (2011), 1335--1338.


\bibitem{jamc2013}
P.~Amodio, T.~Levitina, G.~Settanni and E.B.~Weinm\"uller,
\emph{On the Calculation of the Finite Hankel Transform Eigenfunctions},
J.\  Appl.\  Math.\  Comput.\  \textbf{43} (2013), pp. 151--173,  doi: 10.1007/s12190-013-0657-1.

%\bibitem{aip2013}
%P. ~Amodio, T.~ Levitina, G.~ Settanni, E.B.~ Weinm\"uller,\emph{ Calculations of the Morphology Dependent Resonances}, AIP Conference Proceedings \textbf{1558} (2013), Numerical analysis and applied mathematics - ICNAAM 2013. T.E. Simos, G. Psihoyios and Ch. Tsitouras editors, pp. 750-753, doi: 10.1063/1.4825602

\bibitem{cpc2013}
P.~ Amodio, T.~Levitina, G.~Settanni, E.B.~Weinm\"uller, \emph{Numerical simulation of the whispering gallery modes in prolate spheroids}, accepted to a Computer Physics Communications.

\bibitem{jcam2005}
P.~Amodio and I.~Sgura, %
\emph{High-order finite difference schemes for the solution of
second-order BVPs}, %
J.\ Comput.\ Appl.\ Math.\ \textbf {176} (2005), 59--76.

\bibitem{bit2007}
P.~Amodio and I.~Sgura, %
\emph{High order generalized upwind schemes and numerical
solution of singular perturbation problems}, %
BIT \textbf{47} (2007), 241--257.


\bibitem{AnPa}
A.L.~Andrew and J.W.~Paine, \emph{Correction of Numerov's eigenvalue estimates},
Numer.\ Math.\ \textbf{47} (1985), 289--300.


\bibitem{Ba_th}
P.B.~Bailey, \emph{SLEIGN: an eigenfunction--eigenvalue code for Sturm--Liouville problems},
SAND77-2044, Sandia Laboratories, Albuquerque (1978).

\bibitem{BaGaKaZe}
P.B.~Bailey, B.S.~Garbow, H.G.~Kaper and A.~Zettl,
\emph{Eigenvalue and eigenfunction computations for Sturm-Liouville Problems},
ACM Trans.\ Math.\ Software \textbf{17} (1991), 491--499.

\bibitem{BaEvZe}
P.B.~Bailey, W.N.~Everitt and A.~Zettl, \emph{Algorithm 810: the SLEIGN2 Sturm-Liouville Code},
ACM Trans.\ Math.\ Software \textbf{27} (2001), 143--192.

\bibitem{BaGoSh}
P.B.~Bailey, M.K.~Gordon and L.F.~Shampine,
\emph{Automatic solution of the Sturm-Liouville problem},
ACM Trans.\ Math.\ Software \textbf{4} (1978), 193--208.


\bibitem{BrTr.book}
L.~Brugnano and D.~Trigiante,
\emph{Solving Differential Problems by Multistep Initial and Boundary Value Methods},
Gordon and Breach Science Publishers, Amsterdam, 1998.

\bibitem{AmHiPe}
W.N.~Everitt,
\emph{A catalogue of Sturm-Liouville differential equations},
in: Sturm-Liouville Theory. Past and Present, W.O.~Amrein, A.M.~Hinz and D.P.~Pearson (eds.), Birk\"auser, 2005.

\bibitem{IxDmVb}
L.Gr.~Ixaru, H.~De Meyer and G.~Vanden Berghe, \emph{SLCPM12 -- A program for solving regular Sturm-Liouville problems},
Comput.\ Phys.\ Comm.\ \textbf{118} (1999), 259--277.


\bibitem{LeDaBe}
V.~Ledoux,  M.~Van Daele and G.~Vanden Berghe,
\emph{MATSLISE: a MATLAB package for the numerical solution of Sturm-Liouville and Schr\"odinger equations},
ACM Trans.\ Math.\ Software \textbf{31} (2005), 532--554.


\bibitem{Mar91}
M.~Marletta,
\emph{Certification of algorithm 700: numerical tests of the SLEIGN software for Sturm-Liouville problems}, %
ACM Trans.\ Math.\ Software \textbf{17} (4) (1991), 481--490.

\bibitem{Mar95}
M.~Marletta and J.D.~Pryce, %
\emph{LCNO Sturm-Liouville problems: computational difficulties and examples}, %
Numer.\ Math.\ \textbf{69} (3) (1995), 303--320

\bibitem{PaDeAn}
J.W.~Paine, F.R.~De Hoog and R.S.~Anderssen,
\emph{On the correction of finite difference eigenvalue approxiamtions for Sturm-Liouville problems},
Computing \textbf{26} (1981), 123--139.

\bibitem{Pr}
S.~Pruess,
\emph{Estimating the eigenvalues of Sturm-Liouville problems by approximating the differential equation},
SIAM J.\ Numer.\  Anal.\  \textbf{10} (1973), 55--68.


\bibitem{PrFu}
S.~Pruess and C.T.~Fulton,
\emph{Mathematical software for Sturm-Liouville problems},
ACM Trans.\ Math.\ Software \textbf{19} (1993), 360--376.

\bibitem{PrMa}
S.~Pruess and M.~Marletta,
\emph{Atomatic solution of Sturm-Liouville problems using the Pruess method},
J.\ Comput.\ Appl.\ Math.\ \textbf{39} (1992), 57--78.


\bibitem{Pry}
J.D.~Pryce,  \emph{A test package for Sturm-Liouville solvers},
ACM Trans.\ Math.\ Software \textbf{25} (1) (1999), 21--57.

\bibitem{Pry2}
J.D.~Pryce, \emph{Numerical solution of Sturm-Liouville problems}. Monographs
on Numerical Analysis. Oxford Science Publications. The Clarendon Press,
Oxford University Press, New York, 1993.

\end{thebibliography}
\end{document}